\documentclass{commat}

\usepackage{graphicx}

\DeclareMathAlphabet{\pazocal}{OMS}{zplm}{m}{n}

\title{%
    Exponential stability for a nonlinear porous-elastic system with delay
    }

\author{%
    M.J. Dos Santos, C.A. Raposo, L.G.R. Miranda and B. Feng
    }

\affiliation{
    \address{Manoel Jeremias Dos Santos --
Faculty of Exact Sciences and Technology,  Federal University of Par\'a, 
Abaetetuba, Par\'a,  68440-000, Brazil
        }
    \email{%
    jeremias@ufpa.br
    }
    \address{Carlos Alberto Raposo da Cunha --
    Department of Mathematics, Federal University of Bahia, Salvador, Bahia,  40.170-110, Brazil
        }
    \email{%
    carlos.raposo@ufba.br
    }
    \address{Luiz Gutemberg Ros\'ario Miranda --
    Faculty of Mathematics, Federal University of Par\'a,
Salin\'opolis, Par\'a,  68721-000, Brazil
        }
    \email{%
    luizmiranda@ufpa.br
    }
    \address{Baowei Feng --
    Department of Economic Mathematics, Southwestern University of Finance and Economics,
Chengdu, 611130, People’s Republic of China
        }
    \email{%
    bwfeng@swufe.edu.cn
    }
    }

\abstract{%
    In this work, we consider the existence of a global solution and the exponential decay of a nonlinear porous elastic system with time delay. The nonlinear term, as well as the delay  acting  in the equation of the volume fraction. In order to obtain the existence and uniqueness of a global solution, we will use the semigroup theory of linear operators and under a certain relation involving the coefficients of the system together with a Lyapunov functional, we will establish the exponential decay of the energy associated to the system.
    }



\VOLUME{31}
\YEAR{2023}
\NUMBER{1}
\firstpage{359}
\DOI{https://doi.org/10.46298/cm.10439}

\begin{paper} 

\section{Introduction}

This paper is concerned with a nonlinear porous elastic system with a time delay given by
\begin{align}
\begin{array}{rcl}
    \rho u_{tt}(x,t)\!-\!\mu u_{xx}(x,t)\!-\!b\phi_x(x,t)\!=\!0,\\
    \!\!\!J\phi_{tt}(x,t)\!-\!\delta\phi_{xx}(x,t) \!+\! bu_x(x,t)\!+\!\xi\phi(x,t)\!+\!\mu_1\phi_t(x,t)\!+\!
    \mu_2\phi_t(x,t-\tau)\!+\!f(\phi(x,t))\!=\!0,
\end{array}\label{1:1}
\end{align}
with $(x,t)\in(0,1)\times(0,+\infty)$, were $\rho$, $\mu$, $J$, $\delta$, $\xi$, $\mu_1$, $\mu_2$ are positive constants and $b$ is a non-zero constant satisfying $b^2\leq\mu\xi$.

We consider Dirichlet boundary conditions
\begin{eqnarray}\label{1:2}
u(0,t)=u(1,t)=\phi(0,t)=\phi(1,t) =0 ,\quad t>0
\end{eqnarray}
 and initial conditions
\begin{eqnarray}\label{1:3}
\begin{array}{rl}
u(x,0)=u_0(x),\,\, u_t(x,0)=u_1(x),&  x\in(0,1), \\
\phi(x,0)=\phi_0(x), \,\, \phi_t(x,0)=\phi_0(x),&x\in(0,1).
\end{array}
\end{eqnarray}
Since there exists a delay term in \eqref{1:1} we must provide information for $\phi_t$ in time interval $(0,\tau)$ i.e., in addition to the initial conditions given earlier, we still have
$$
\phi_t(x,t-\tau)=f_0(x,t-\tau), \,(x,t)\in(0,1)\times(0,\tau).
$$
where $(u_0,u_1,\phi_0,\phi_1,f_0)$ will be taken in an appropriate space.

The importance of porous elastic materials lies in the fact that they are present in various fields of human activity such as the petroleum industry, material science, soil mechanics, foundation engineering, powder technology, biology and others (cf. \cite{IESANBOOK}, \cite{STRAUGHAN}).

Among the various theories dealing with a porous material, we can find a linear theory proposed by Cowin and Nunziato \cite{COWIN-NUNZIATO}, \cite{NUNZIATO-COWIN} which is a generalization of the elastic theory for materials with voids, considering besides the material elasticity property, the volume fraction of the voids in the material. In this theory, the bulk density $\rho=\rho(x,t)$ is given by the product of matrix density of the material $\gamma=\gamma(x,t)$ and the volume fraction $\nu=\nu(x,t)$
$$\rho(x,t)=\gamma(x,t)\nu(x,t).$$
They also consider a reference configuration (generally considered  as  an initial configuration)
$$\rho_0(x)=\gamma_0(x)\nu_0(x).$$

Let $u_i=u_i(x,t)$ denote the components of the displacement vector field and so the components of the infinitesimal strain field are given by
$$e_{ij}=\frac{1}{2}(u_{i,j}+u_{j,i}),$$
where the comma after the letter indicates the partial derivative with respect to the indicated coordinate. In addition, $\phi=\phi(x,t) $ represents the change in volume fraction with respect to the reference configuration.
In a framework of evolution equations and in a setting of three-dimensional theory, the porous-elastic theory is  described  as
\begin{eqnarray*}
\rho\ddot{u}_i&=&T_{ij,j}+\rho f_i,\\
\rho k\ddot{\phi}&=&h_{i,i}+g+\rho\ell,
\end{eqnarray*}
which are called the balance of linear  momentum  and the balance of  equilibrated force  equations, respectively. In the above system, $T_{ji}$ are the components of the stress tensor, $f_i$ is the body force vector, $h_i$ are the components of the equilibrated stress vector, $k$ is the equilibrated inertia, $g$ is the intrinsic equilibrated body force and $\ell$ is the extrinsic equilibrated body force.

The constitutive equations for homogeneous and isotropic elastic bodies are (cf. \cite{COWIN-NUNZIATO}):
\begin{eqnarray*}
T_{ij}&=&\lambda \delta_{ij} e_{rr}+2\mu e_{ij}+\beta\phi\delta_{ij},\\
h_i&=&\alpha\phi_{,i},\\
g&=&-\omega\dot{\phi}-\xi\phi-\beta e_{rr},
\end{eqnarray*}
where $\lambda$, $\mu$, $\beta$, $\alpha$, $\omega$, $\xi$ and $\omega$, are constitutive constants that depend on the reference state $\nu_0$ and $\delta_{ij}$ is Kronecker's delta. The necessary and sufficient conditions for the internal energy density to be positive definite quadratic form are (cf.\cite{COWIN-NUNZIATO})
$$\mu>0, \ \alpha>0, \ \xi>0, \ \kappa>0, \ \omega>0, \ 3\lambda+2\mu>0 \ \mbox{and} \ 3b^2\leq(3\lambda+2\mu)\xi.$$

\subsection{The stabilization scenario}

The stabilization study constitutes an important issue in research areas into mathematics and engineering and has been aiming for investigation for years. Naturally, significant mathematical properties are extracted from undamped partial
differential equations, and we pointed out the important mathematical properties from wave equations, plate equations, shell equations and beam (plane and curved) beams \cite{GRAFFBOOK}, \cite{LAGNESE-LIONS-1}. However, dissipative mechanisms make more realistic the phenomenons translated in terms of partial differential equations. See for example the classical books due to Komornik \cite{KOMORNIK} and Lions and Lagnese  \cite{LAGNESEBOOK}, \cite{LIONSBOOK} on boundary stabilization.  In this direction, we will say that the energy of the solutions of dissipative partial differential equations
decay exponentially if they are controlled by an exponential negative. On the contrary, we will say that the energy of the solution
slowly decays. Thus, in both cases is  important to determine conditions to achieve some type of decay.

To the best knowledge, the analysis of the temporal decay in one-dimensional porous-elasticity theory was first studied by
Quintanilla \cite{QUINTANILLA}. He
showed, by taking into account a particular set of solutions and by using the Routh-Hurwitz theorem, that the porous viscosity acting on the second equation of the following porous elastic system
\begin{eqnarray*}
 \rho_0 u_{tt}-\mu u_{xx}-\beta\varphi_x&=&0,\\
 \rho_0\kappa \varphi_{tt}-\alpha \varphi_{xx}+\beta u_x+\xi \varphi+\tau \varphi_t&=&0,
\end{eqnarray*}
was not powerful enough to obtain the exponential decay of solutions when $\alpha\neq\mu k$.
For this problem was used the denomination of slow decay was to characterize the decay obtained.

Damping effects like temperature, elastic viscosity, porous viscosity and micro-tem\-per\-a\-ture acting on porous-elasticity equations were considered in the past years and obtained several contributions to porous-elasticity systems with an emphasis on lack of exponential decay, the exponential decay, the polynomial decay, the optimality of the polynomial decay and  numerical computations results, see
\cite{Casas_PS}, \cite{CASAS-QUINTANILLA-2005-2}, \cite{MAGANA-QUINTANILLA-2006}, \cite{MAGANA-QUINTANILLA-2007}, 
\cite{PAMPLONA-RIVERA-QUINTANILLA-2009}, \cite{PAMPLONA-RIVERA-QUINTANILLA-2011},  \cite{SANTOS-ALMEIDAJR-2016}, \cite{SANTOS-ALMEIDAJR-2017}, \cite{SANTOS-CAMPELO-ALMEIDAJR-2016}, \cite{SANTOS-CAMPELO-ALMEIDAJR-2017}.

A porous thermoelastic system was investigated in \cite{FAREH-MESSAOUDI-2017} where the heat conduction is given by Cattaneo's law and was proved an exponential decay result. Analyticity of the semigroup associated with the transmission problem in porous elasticity and delay was studied in \cite{RAPOSO-APALARA-RIBEIRO}. It is interesting to note a particular case of the equations \eqref{1:1}, when $k:=\mu=\xi=b$ and $\rho_0:=\rho$ e $\rho_1:=\rho k$, we obtain the following nonlinear Timoshenko system that modeling a thin elastic beam with delay
\begin{eqnarray}
\begin{array}{rcl}
    \rho_0 u_{tt}-k(u_{x}+\phi)_x&=&0,\\
    \rho_1\phi_{tt}-\delta\phi_{xx}+k(u_x+\phi)+\mu_1\phi_t+\mu_2\phi_t(x,t-\tau)+f(\phi)&=&0.
\end{array}\label{1:4}
\end{eqnarray}
The exponential decay for the system \eqref{1:4} was studied in \cite{FENG-PELICER} which extended the result previously obtained in \cite{HOUARI-LASKRI}.  In both cases, linear and non-linear, exponential decay was obtained through equality
\begin{eqnarray}
\frac{\rho_0}{k}=\frac{\rho_1}{\delta}.\label{1:5}
\end{eqnarray}
This paper intends to answer the following question: Is it possible to add terms of delay, friction and forcing in the porous elastic system similar to that found in \eqref{1:4} and also to obtain the exponential decay? The answer is yes, and in this situation, a relation between the coefficients will appear similar to \eqref{1:5}, in the porous elastic case \eqref{1:1} the equality
\begin{eqnarray}
\frac{\rho}{\mu}=\frac{J}{\delta}\label{1:6}
\end{eqnarray}
will play an important role.

This article is divided as follows: in section 2 we will establish the notations, as well as put the system \eqref{1:1}-\eqref{1:3} in the form of an abstract Cauchy problem. In section 3 we will establish the existence of a global solution for the system \eqref{1:1}-\eqref{1:1} and in section 4 we will deal with the main result of this work, which is the exponential decay of the energy of the solution for the problem \eqref{1:1}-\eqref{1:3}.

\section{Preliminaries and main results}
In this work, we will use the standard notations $L^p(0,1)$ with $1\leq p\leq\infty$ for the Lebesgue spaces  and $H_0^1(0,1)$ for Sobolev space. The norm in $L^p(0,1)$ will be represented by $\|\cdot\|_p$ and we will write $\|\cdot\|$ instead of $\|\cdot\|_2$ to norm in $L^2(0,1)$.

For the forcing term $f(\phi)$, we assume $f:\mathbb{R}\to\mathbb{R}$ and that there are constants $k_0>0$ and $\theta>0$ satisfying
\begin{eqnarray}
|f(x)-f(y)|\leq k_0(|x|^\theta+|y|^\theta)|x-y|\quad \forall x,y\in\mathbb{R}. \label{2:1}
\end{eqnarray}
In addition, we will also assume that
\begin{eqnarray}
0\leq \hat{f}(x)\leq f(x)x\quad\forall x\in\mathbb{R},\label{2:2}
\end{eqnarray}
with $\hat{f}(x)=\displaystyle\int_0^xf(s)ds$.

Using an argument found in \cite{FENG-PELICER}, \cite{HOUARI-LASKRI}, \cite{NICAISE-PIGNOTTI-2006}, \cite{NICAISE-PIGNOTTI-2008}, in order to deal with the delay term, we define the following variable:
\begin{eqnarray}
z(x,y,t)=\phi_t(x,t-\tau y),\quad x\in(0,1), \ y\in(0,1), \ t>0.\label{2:3}
\end{eqnarray}
Then it is easy to verify
\begin{eqnarray}
\tau z_t(x,y,t)+z_y(x,y,t)=0,\quad \forall (x,y,t)\in(0,1)\times(0,1)\times(0,\infty).\label{2:4}
\end{eqnarray}
Thus, the system \eqref{1:1}-\eqref{1:3} takes the following form
\begin{align}
    \rho u_{tt}(x,t)-\mu u_{xx}(x,t)-b\phi_x(x,t)\!=\!0,\label{2:5}\\
    \!\!\!\!J\phi_{tt}(x,t)-\delta\phi_{xx}(x,t)+bu_x(x,t)+\xi\phi(x,t)+\mu_1\phi_t(x,t)+
    \mu_2z(x,1,t)+f(\phi(x,t))\!=\!0,\label{2:6}\\
\tau z_t(x,y,t)+z_y(x,y,t)\!=\!0,\label{2:7}
\end{align}
with $x\in(0,1)$, $y\in(0,1)$ and $t>0$. In addition, boundary conditions are given by
\begin{eqnarray}
\begin{array}{rl}
u(0,t)=u(1,t)=\phi(0,t)=\phi(1,t)=0, &t>0,\\
z(x,0,t)=\phi_t(x,t),& (x,t)\in (0,1)\times(0,+\infty)
\end{array}\label{2:8}
\end{eqnarray}
and initial conditions
\begin{align}
\begin{array}{rl}
u(x,0)=u_0(x),\, u_t(x,0)=u_1(x), \, \phi(x,0)=\phi_0(x),\, \phi_t(x,0)=\phi_0(x),x\in(0,1),\\
z(x,1,0)=f_0(x,t-\tau),(x,t)\in(0,1)\times(0,\tau).
\end{array}\label{2:9}
\end{align}

Considering $U(t)=(u(t),u_t(t),\phi(t),\phi_t(t),z(t)$ and $U_0=(u_0,u_1,\phi_0,\phi_1,f_0)$, the system or \eqref{2:5}-\eqref{2:9} is reduced to the following abstract first order evolution problem:
\begin{equation}
\left\{\begin{array}{rcl}
\dot{U}(t)&=&\pazocal{A}U(t)+\pazocal{F}(U(t)), \quad t>0,\\
U(0)&=&U_0,
\end{array}\right.\label{3:1}
\end{equation}
where\,\, $\dot{}=\frac{d}{dt}$,  and
\begin{equation}
\pazocal{A}U=\left(\begin{array}{c}
v\\
\frac{\mu}{\rho}u_{xx}+\frac{b}{\rho}\phi_x\\
\psi\\
\frac{\delta}{J}\phi_{xx}-\frac{b}{J}u_x-\frac{\xi}{J}\phi-\frac{\mu_1}{J}\psi-\frac{\mu_2}{J}z(\cdot,1)\\
-\frac{1}{\tau}z_y
\end{array}\right), \quad
\pazocal{F}(U)=\left(\begin{array}{c}
0\\ 0\\ 0\\ -\frac{1}{J}f(\phi)\\ 0
\end{array}\right)\label{3:1a}
\end{equation}
with the domain
\begin{eqnarray}
D(\pazocal{A})=\{(u,v,\phi,\psi,z)^T\in H; \ \psi(x)=z(x,0), \ x\in(0,1)\},\label{3:2}
\end{eqnarray}
where 
\begin{align}
H=(H^2(0,1)\cap H^1_0(0,1))\!\times\! H_0^1(0,1)\!\times\!(H^2(0,1)\cap H_0^1(0,1))
\!\times\! H_0^1(0,1)\!\times\! L^2(0,1;H^1(0,1)).
\end{align}
We define the energy space $\mathcal{H}$ by
\begin{eqnarray}
\mathcal{H}:=H^1_0(0,1)\times L^2(0,1)\times H_0^1(0,1)\times L^2(0,1)\times L^2(0,1;L^2(0,1)).\label{3:3}
\end{eqnarray}

For $U=(u,v,\phi,\psi,z)^T$, $V=(\tilde{u},\tilde{v},\tilde{\phi},\tilde{\psi},\tilde{z})^T$ in $\mathcal{H}$ and for $\eta$ a positive constant satisfying
\begin{eqnarray}
\tau\mu_2\leq \eta\leq \tau (2\mu_1-\mu_2),\label{3:4}
\end{eqnarray}
we equip $\mathcal{H}$ with the following inner product:
\begin{align}
\langle U,V\rangle_{\mathcal{H}}=&\int_0^1[\rho v\tilde{v}+J\psi\tilde{\psi}+\delta\phi_x\tilde{\phi}_x+\mu u_x\tilde{u}_x+b(u_x\tilde{\phi}+\tilde{u}_x\phi)+\xi\phi\tilde{\phi}]dx \nonumber \\
&\,\,+\eta\int_0^1\int_0^1z(x,y)\tilde{z}(x,y)dy dx\label{3:5}
\end{align}
and so
\begin{eqnarray}
\|U\|^2_{\mathcal{H}}=\int_0^1[\rho v^2+J\psi^2+\delta\phi_x^2+\mu u_x^2+2bu_x\phi+\xi\phi^2]dx
+\eta\int_0^1\int_0^1z^2(x,y)dy dx.\label{3:6}
\end{eqnarray}

\begin{remark}
Since $b^2\leq \mu\xi$ and
\begin{eqnarray}
\mu u_x^2+2bu_x\phi+\xi\phi^2=\bigg(\frac{b}{\sqrt{\xi}}u_x+\sqrt{\xi}\phi\bigg)^	2+\bigg(\mu-\frac{b^2}{\xi}\bigg)
u_x^2,
\label{3:18}
\end{eqnarray}
it is easy to see that $\|U\|^2_{\mathcal{H}}\geq 0$ for all $U\in\mathcal{H}$.
\end{remark}

\begin{definition}
Let $S(t)$ be a  $C_0$-semigroup of contractions on $\mathcal{H}$. A continuous function $U\,:\, (0,T] \rightarrow \mathcal{H}$ is a mild solution to the problem \eqref{3:1} if it satisfies
\begin{eqnarray}
U(t)=S(t)U_0+\int_0^t S(t-s)f(U(s))ds,\,\,\,0 < t \leq T,
\end{eqnarray}
where $f$ is locally Lipschitz on $\mathcal{H}$.
\end{definition}

The first main result regarding the existence and uniqueness of the global solution to the system \eqref{3:1} and which will play an important role in the result of the exponential decay is

\begin{theorem}\label{T:3:1}
Assume that \eqref{2:1}-\eqref{2:2} and $\mu_2\leq\mu_1$ hold, then we have the following results.
\begin{description}
  \item[(i)] If $U_0\in\mathcal{H}$, then problem \eqref{3:1} has a unique mild solution $U\in C([0,\infty),\mathcal{H})$ with $U(0)=U_0$;
  \item[(ii)] If $U_1$ and $U_2$ are two mild solutions of problem \eqref{3:1}, then there exists a positive constant $C_0=C(U_1(0),U_2(0))$ such that
  \begin{eqnarray}
  \|U_1(t)-U_2(t)\|_{\mathcal{H}}\leq e^{C_0T}\|U_1(0)-U_2(0)\|_{\mathcal{H}}, \quad \forall 0\leq t\leq T.\label{3:24}
  \end{eqnarray}
  \item[(iii)] If $U_0\in D(\pazocal{A})$, then the above mild solution is a strong solution.
\end{description}
\end{theorem}

With the solution $U(t)=(u(t),u_t(t),\phi(t),\phi_t(t),z(t))$  of \eqref{2:5}-\eqref{2:9} obtained by Theorem~\ref{T:3:1}, we can define the energy associated with this solution as follows:

\begin{equation}
E(t)=\frac{1}{2}\int_0^1[\rho u_t^2+J\phi_t^2+\delta\phi_x^2+\mu u_x^2+2bu_x\phi+\xi\phi^2+2\hat{f}(\phi)]dx\\+
\frac{\eta}{2}\int_0^1\int_0^1z^2(x,y,t)dydx.
\label{3:17}
\end{equation}

The second main result of this work, tells us about the energy behavior of the solution over time and its statement is

\begin{theorem}\label{T:4:1}
Assume that \eqref{2:1}-\eqref{2:2} and $\mu_2<\mu_1$ hold. Assume that \eqref{1:6} also holds. Then, with respect to mild solution, there exist $C>0$ and $\gamma>0$ such that
\begin{eqnarray}
E(t)\leq Ce^{-\gamma t}, \quad t\geq 0.\label{4:39}
\end{eqnarray}
\end{theorem}

\section{Well-posedness of the problem}

\begin{lemma}\label{L:3:1}
The operator $\pazocal{A}$ defined in \eqref{3:1a} is the infinitesimal generator of a $C_0$-sem\-i\-group of contractions on $\mathcal{H}$.
\end{lemma}
\begin{proof}
Firstly it is not difficult to see that $D(\pazocal{A})$ is dense in $\mathcal{H}$. In addition, it follows of definition of $\pazocal{A}$, \eqref{3:5} and by integration by parts that for any $U=(u,v,\phi,\psi,z)\in D(\pazocal{A})$
\begin{eqnarray}
\langle \pazocal{A}U,U\rangle_{\mathcal{H}}= \bigg(-\mu_1+\frac{\eta}{2\tau}\bigg)\int_0^1\psi^2dx-\mu_2\int_0^1z(x,1)\psi dx-\frac{\eta}{2\tau}\int_0^1z^2(x,1)dx,
\end{eqnarray}
and applying H\"older's inequality we have
\begin{eqnarray}
\langle \pazocal{A}U,U\rangle_{\mathcal{H}}\leq\bigg(-\mu_1+\frac{\eta}{2\tau}+\frac{\mu_2}{2}\bigg)\int_0^1\psi^2dx
+\bigg(-\frac{\eta}{2\tau}+\frac{\mu_2}{2}\bigg)\int_0^1z^2(x,1)dx.
\end{eqnarray}
Keeping in mind the condition \eqref{3:4}, we observe that
\begin{eqnarray}
-\mu_1+\frac{\eta}{2\tau}+\frac{\mu_2}{2}\leq 0\quad\text{and}\quad -\frac{\eta}{2\tau}+\frac{\mu_2}{2}\leq0,
\end{eqnarray}
with implies that the operator $\pazocal{A}$ is dissipative.

We will now prove that $0\in\varrho(\pazocal{A})$, where $\varrho(\pazocal{A})$ is resolvent set of $\pazocal{A}$, that is, given $V=(f_1,f_2,f_3,f_4,f_5)\in\mathcal{H}$, we must find $U=(u,v,\phi,\psi,z)\in D(\pazocal{A})$ solution of
\begin{eqnarray}
-\pazocal{A}U=V, \quad\text{with}\quad \|U\|_{\mathcal{H}}\leq C\|V\|_{\mathcal{H}},
\end{eqnarray}
where $C$ is a constant that is independent of $U$ and $V$. The previous equality results in the following system:
\begin{eqnarray}
-v&=&f_1, \label{3:8}\\
-\mu u_{xx}-b\phi_x&=&\rho f_2, \label{3:9}\\
-\psi&=&f_3, \label{3:10}\\
-\delta\phi_{xx}+bu_x+\xi\phi+\mu_1\psi+\mu_2z(\cdot,1)&=&Jf_4, \label{3:11}\\
z_y&=&-\tau f_5.\label{3:12}
\end{eqnarray}
By \eqref{3:8} and \eqref{3:10} we have, $v\in H_0^1(0,1)$ and $\psi\in H_0^1(0,1)$ and by \eqref{3:12}
\begin{eqnarray}
z(\cdot,y)&=&-\tau\int_0^y f_5(\cdot,s)ds-f_3(\cdot)\in L^2(0,1;H^1(0,1)).\label{3:13}
\end{eqnarray}
Substituting \eqref{3:10} and \eqref{3:12} into \eqref{3:11} we obtain the system
\begin{eqnarray}
\begin{array}{rcl}
-\mu u_{xx}-b\phi_x&=&F_1\in L^2(0,1),\\
-\delta\phi_{xx}+bu_x+\xi\phi&=&F_2\in L^2(0,1),
\end{array}\label{3:14}
\end{eqnarray}
where
\begin{eqnarray*}
F_1&=&\rho f_2,\\
F_2&=&(\mu_1+1)f_3+\mu_2\tau\int_0^1f_5(\cdot,s)ds+Jf_4,
\end{eqnarray*}
with boundary conditions
\begin{eqnarray}
u(0)=u(1)=\phi(0)=\phi(1)=0.\label{3:15}
\end{eqnarray}
Solving system \eqref{3:14}-\eqref{3:15} is equivalent to finding $(u,\phi)\in H_0^1(0,1)\times H_0^1(0,1)$ such that
\begin{eqnarray}
\int_0^1[\mu u_x\bar{u}_x+b(u_x\bar{\phi}+\bar{u}_x\phi)+\xi\phi\bar{\phi}+
\delta\phi_x\bar{\phi}_x]dx=\int_0^1F_1\bar{u}dx+\int_0^1F_2\bar{\phi}dx\label{3:16},
\end{eqnarray}
for all $(\bar{u},\bar{\phi})\in H_0^1(0,1)\times H_0^1(0,1)$. Problem \eqref{3:16} is equivalent to
\begin{eqnarray}
a((u,\phi),(\bar{u},\bar{\phi}))=h(\bar{u},\bar{\phi})
\end{eqnarray}
where the bilinear form $a:(H_0^1(0,1)\times H_0^1(0,1))^2\to\mathbb{R}$ and the linear form $$h:H_0^1(0,1)\times H_0^1(0,1)\to\mathbb{R}$$ are given by
\begin{eqnarray}
a((u,\phi),(\bar{u},\bar{\phi}))=\int_0^1[\mu u_x\bar{u}_x+b(u_x\bar{\phi}+\bar{u}_x\phi)+\xi\phi\bar{\phi}+
\delta\phi_x\bar{\phi}_x]dx
\end{eqnarray}
and
\begin{eqnarray}
h(\bar{u},\bar{\phi})=\int_0^1F_1\bar{u}dx+\int_0^1F_2\bar{\phi}dx.
\end{eqnarray}
The bilinear form $a$ is coercive because thanks to Young's inequality we have
$$\begin{array}{lll}
a((u,\phi),(u,\phi))&=&\int_0^1[\mu u_x^2+2bu_x\phi+\xi\phi^2+\delta\phi_x^2]dx\\
&\geq&\mu\int_0^1u_x^2dx-\frac{b^2}{\xi}\int_0^1u_x^2dx\\
& &-\xi\int_0^1\phi^2dx+\xi\int_0^1\phi^2dx+\delta\int_0^1\phi_x^2dx\\
&\geq&\bigg(\frac{\mu\xi-b^2}{\xi}\bigg)\int_0^1u_x^2dx+\delta\int_0^1\phi_x^2dx\\
&\geq&C_1(\|u_x\|^2+\|\phi_x\|^2),
\end{array}$$
where $C_1=\min\{(\mu\xi-b^2)/\xi,\delta\}$. It is not difficult to verify that $a$ and $h$ are continuous. So applying the Lax-Milgram's Theorem, we deduce that problem \eqref{3:14} admits a unique solution $(u,\phi)\in H_0^1(0,1)\times H_0^1(0,1)$. Thanks to the classical elliptic regularity (c.f. \cite{NIRENBERG}) we obtain $(u,\phi)\in H^2(0,1)\times H^2(0,1)$ and $C>0$. Therefore, $0\in\varrho(\pazocal{A})$. Consequently, the result of the lemma follows from Theorem 1.2.4 in \cite{LIU-ZHANG}.
\end{proof}

\begin{lemma}\label{L:3:2}
The operator $\pazocal{F}$ defined in \eqref{3:1} is locally Lipschitz in $\mathcal{H}$.
\end{lemma}
\begin{proof}
Let $U=(u,v,\phi,\psi,z)$ and $V=(\tilde{u},\tilde{v},\tilde{\phi},\tilde{\psi},\tilde{z})$ be in $\mathcal{H}$, with
\begin{eqnarray}
\|U\|_{\mathcal{H}},\|V\|_{\mathcal{H}}\leq M, \quad M>0.
\end{eqnarray}
By using \eqref{2:1}, Hölder's inequality and the embedding $H^1(0,1)\hookrightarrow L^{2(\theta+1)}(0,1)$ we obtain
\begin{align*}
\|F(U)-F(V)\|^2_{\mathcal{H}} 
&= J^{-1}\int_0^1|f(\phi)-f(\tilde{\phi})|^2 \\
&\leq \ J^{-1} k_0^2\int_0^1(|\phi|^\theta+|\tilde{\phi}|^\theta )^2|\phi-\tilde{\phi}|^2dx\\
&\leq C_0\bigg(\int_0^1(|\phi|^\theta+|\tilde{\phi}|^\theta )^\frac{2(\theta+1)}{\theta}\bigg)^{\frac{\theta}{\theta+1}}\bigg(\int_0^1|\phi-\tilde{\phi}|^{2(\theta+1)}dx\bigg)^{\frac{1}{\theta+1}}\\
&\leq C_0(\|\phi\|_{2(\theta+1)}^{2\theta}+\|\tilde{\phi}\|_{2(\theta+1)}^{2\theta})\|\phi-\tilde{\phi}\|_{2(\theta+1)}^2\\
&\leq C_0(\|\phi\|^{2\theta}+\|\tilde{\phi}\|^{2\theta})\|\phi-\tilde{\phi}\|^2\\
&\leq C_0(M)\|U-V\|_{\pazocal{H}}^2,
\end{align*}
where $C_0$ is a positive constant depending on $M$. 
\end{proof}

Based on the energy defined in \eqref{3:17}, we have the following result

\begin{lemma}\label{L:3:3}
The functional energy $E(t)$ defined in \eqref{3:17} of a solution $$U(t)=(u(t),u_t(t),\phi(t),\phi_t(t),z(t))$$ of \eqref{2:5}-\eqref{2:9} is  non-increasing function. More precisely, there exists a constant $C_E>0$ such that for all $t\geq 0$
\begin{eqnarray}
\frac{dE}{dt}(t)\leq -C_E\bigg(\int_0^1\phi_t^2dx+\int_0^1z^2(x,1,t)dx\bigg)\leq 0.\label{3:20}
\end{eqnarray}
\end{lemma}
\begin{proof}
Multiplying \eqref{2:5} by $u_t$, \eqref{2:6} by $\phi_t$, integrating the result over $[0,1]$ with respect to $x$ and using Young's inequality, we obtain
\begin{align}
\begin{split}
\frac{1}{2}\frac{d}{dt}\bigg(\int_0^1[\rho u_t^2+J\phi_t^2+\delta\phi_x^2&+\mu u_x^2+2bu_x\phi+\xi\phi^2 
+2\hat{f}(\phi)]dx\bigg)\\
&= -\mu_1\int_0^1 \phi_t^2dx-\mu_2\int_0^1z(x,1,t)\phi_tdx\\
&\leq-\bigg(\mu_1-\frac{\mu_2}{2}\bigg)\int_0^1\phi_t^2dx+\frac{\mu_2}{2}\int_0^1z^2(x,1,t)dx.
\end{split}\label{3:21}
\end{align}
Multiplying \eqref{2:7} by $\dfrac{\eta}{\tau}z$ and integrating the result over $[0,1]\times[0,1]$ with respect $y$ and $x$, respectively we have
\begin{align}
\begin{split}
\frac{\eta}{2}\frac{d}{dt}\int_0^1\int_0^1z^2(x,y,t)dy dx&=-\frac{\eta}{2\tau}\int_0^1\int_0^1\frac{\partial}{\partial y}z^2(x,y,t)dy dx  \\
&=\frac{\eta}{2\tau}\int_0^1(z^2(x,0,t)-z^2(x,1,t))dx \\
&=\frac{\eta}{2\tau}\int_0^1(\phi_t^2-z^2(x,1,t))dx,  
\end{split}\label{3:22}
\end{align}
adding \eqref{3:21} and \eqref{3:22} we obtain
\begin{eqnarray}
\frac{dE}{dt}(t)\leq-\bigg(\mu_1-\frac{\mu_2}{2}-\frac{\eta}{2\tau}\bigg)\int_0^1\phi_t^2dx-
\bigg(\frac{\eta}{2\tau}-\frac{\mu_2}{2}\bigg)\int_0^1z^2(x,1,t)dx.
\end{eqnarray}
Considering
\begin{eqnarray}
C_E=\min\bigg\{\mu_1-\frac{\mu_2}{2}-\frac{\eta}{2\tau},\frac{\eta}{2\tau}-\frac{\mu_2}{2}\bigg\},
\end{eqnarray}
we obtain \eqref{3:20}.
\end{proof}

\begin{remark}
It follows from Theorem \ref{L:3:3} that
\begin{eqnarray}
E(t)\leq E(0),\quad \forall t\geq 0.\label{3:23-1}
\end{eqnarray}
Furthermore, if $U(t)=(u(t),u_t(t),\phi(t),\phi_t(t),z(t))$ is solution of \eqref{2:5}-\eqref{2:9}, then it follows from \eqref{2:2}, \eqref{3:6} and \eqref{3:17}
\begin{align}
E(t) 
&\geq \frac{1}{2}\bigg\{\int_0^1[\rho u_t^2+J\phi_t^2+\delta\phi_x^2+\mu u_x^2+2bu_x\phi+\xi\phi^2dx+ \eta\int_0^1\int_0^1z^2(x,y,t)dydx\bigg\} \label{3:23-2} \\
&\geq \frac{1}{2}\|U(t)\|_{\mathcal{H}}^2. \nonumber
\end{align}
\end{remark}

\begin{remark}
By Using Young's and Poincaré's inequalities in
\small
\begin{eqnarray*}
\frac{b^2}{\xi}\int_0^1u^2_xdx=\int_0^1\bigg(\frac{b}{\sqrt{\xi}}u_x+\sqrt{\xi}\phi\bigg)^2dx-\int_0^1bu_x\phi dx
-\sqrt{\xi}\int_0^1\phi\bigg(\frac{b}{\sqrt{\xi}}u_x+
\sqrt{\xi}\phi\bigg)dx
\end{eqnarray*}\normalsize
we obtain,
\begin{eqnarray}
\int_0^1u^2_xdx\leq\frac{3\xi}{b^2}\int_0^1\bigg(\frac{b}{\sqrt{\xi}}u_x+\sqrt{\xi}\phi\bigg)^2dx+
\frac{2\xi^2c_p}{b^2}\int_0^1\phi_x^2dx.\label{3:19}
\end{eqnarray}
\end{remark}

We are now able to prove the main result of this section which is:

\begin{proof}[\textit{Proof of Theorem \ref{T:3:1}}]
\noindent\textbf{i)} Since $\pazocal{A}$ generates a $C_0$-semigroup of contractions (Lemma \ref{L:3:1}) and $\pazocal{F}$ is locally Lipschitz (Lemma \ref{L:3:2}) follows from Pazy \cite{PAZY}, Theorem 6.1.4 that \eqref{3:1} has a unique mild solution
\begin{eqnarray}
U(t)=e^{\pazocal{A}t}U_0+\int_0^te^{\pazocal{A}(t-s)}\pazocal{F}(U(s))ds,
\end{eqnarray}
defined in $[0,t_{max})$, where $t_{max}$ depend on $U_0$. Moreover, if $t_{max}<\infty$, then
\begin{eqnarray}
\lim_{t\to t_{\text{max}}}\|U(t)\|_{\mathcal{H}}=+\infty.
\end{eqnarray}
However, from \eqref{3:23-1} and \eqref{3:23-2} we have
\begin{eqnarray}
\|U(t)\|_{\mathcal{H}}^2\leq 2E(0)\quad\text{for all}\quad t\geq 0.
\end{eqnarray}
Therefore, $t_{max}=+\infty$.

\noindent\textbf{ii)} It is not difficult to obtain \eqref{3:24}, just consider a standard procedure found in Pazy \cite{PAZY}, Theorem 6.1.2, p. 184, which consists in using the locally Lipschitz condition of $\pazocal{F}$ and the  Gr\"{o}nwall's Inequality.

\noindent\textbf{iii)} By using Theorem 6.1.5 in Pazy \cite{PAZY}, we know that any mild solutions with initial data in $D(\pazocal{A})$ are strong solution.
\end{proof}

\section{Exponential stability}

In this section, we will prove the second main result of this work, Theorem \ref{T:4:1}, which will be divided into the following lemmas.

\begin{lemma}\label{L:4:1}
Let $U(t)=(u(t),u_t(t),\phi(t),\phi_t(t),z(t))$ be the solution to \eqref{2:5}-\eqref{2:9}. Then there is $\delta_0>0$ (depending only on the coefficients of \eqref{2:5}-\eqref{2:7}) such that
\begin{equation}
\begin{split}
E(t)\!\leq\! \delta_0\bigg\{\!\!\int_0^1\!\!\bigg[u_t^2+\phi_t^2+\phi_x^2+
\bigg(\frac{b}{\sqrt{\xi}}u_x+\sqrt{\xi}\phi\bigg)^2\!\!+
\hat{f}(\phi)\bigg]dx+\!\int_0^1\!\!\int_0^1 \!\!\! z^2(x,y,t)dy dx\bigg\}.
\end{split}\label{4:1}
\end{equation}
\end{lemma}
\begin{proof}
It is enough to consider \eqref{3:17}, \eqref{3:18} and \eqref{3:19}.
\end{proof}

Let us now consider a solution $U(t)=(u(t),u_t(t),\phi(t),\phi_t(t),z(t))$ to \eqref{2:5}-\eqref{2:9} and for this solution we define the following functionals:
\begin{eqnarray}
I_1(t)&:=&-\int_0^1(\rho uu_t+J\phi\phi_t)dx-\frac{\mu_1}{2}\int_0^1\phi^2dx,\label{4:2}\\
I_2(t)&:=&\int_0^1(J\phi_t\phi+\rho u_tw)dx+\frac{\mu_1}{2}\int_0^1\phi^2dx,\label{4:5}
\end{eqnarray}
where $w$ is the solution of
\begin{eqnarray}
-w_x=\frac{b}{\mu}\phi\quad w(0)=w(1)=0,\label{4:6}
\end{eqnarray}
\begin{eqnarray}
I_3(t)&:=&J\int_0^1\bigg(\frac{b}{\sqrt{\xi}}u_x+\sqrt{\xi}\phi\bigg)\phi_tdx+\frac{b J}{\sqrt{\xi}}\int_0^1\phi_x u_tdx,\\
I_4(t)&:=&\int_0^1\int_0^1e^{-2\tau y}z^2(x,y,t)dydx.\label{4:29}
\end{eqnarray}
Thereby, we have the following results.
\begin{lemma}\label{L:4:2}
For any $\varepsilon>0$ the functional $I_1$ satisfies,
\begin{align}
\frac{d}{dt}I_1(t)
\leq{ }&{ }-\int_0^1(\rho u_t^2+J\phi_t^2)dx +\bigg[\delta+\varepsilon+c_1+\frac{2\xi^2c_p}{b^2}\bigg(\mu-\frac{b^2}{\xi}\bigg)\bigg]\int_0^1\phi_x^2dx \label{4:3} \\
&{ }+\bigg[3\bigg(\frac{\mu\xi}{b^2}-1\bigg)+1\bigg]\int_0^1\bigg(\frac{b}{\sqrt{\xi}}u_x+\sqrt{\xi}\phi\bigg)^2dx + \frac{\mu_2^2c_p}{4\varepsilon}\int_0^1z^2(x,1,t)dx, \nonumber
\end{align}
where $c_1>0$ is a constant depending on Poincar\'e's constant $c_p$.
\end{lemma}

\begin{proof}
Taking the derivative of $I_1$, we have
\begin{eqnarray*}
\frac{d}{dt}I_1(t)=-\int_0^1(\rho u^2_t+J\phi_t^2)dx-\int_0^1(\rho u_{tt}u+J\phi_{tt}\phi)dx-\mu_1\int_0^1\phi\phi_tdx.
\end{eqnarray*}
Using \eqref{2:5}, \eqref{2:6} and \eqref{3:18} we arrive at
\begin{equation}
\begin{split}
\frac{d}{dt}I_1(t)=-\int_0^1(\rho u^2_t+J\phi_t^2)dx+\delta\int_0^1\phi_x^2dx+\int_0^1\bigg(\frac{b}{\sqrt{\xi}}u_x+\sqrt{\xi}\phi\bigg)^2dx\\
+\bigg(\mu-\frac{b^2}{\xi}\bigg)\int_0^1u_x^2dx+\mu_2\int_0^1z(x,1,t)\phi dx+\int_0^1f(\phi)\phi dx.\label{4:4}
\end{split}
\end{equation}
It follows from Young's and Poincaré's inequalities that for any $\varepsilon>0$,
\begin{align}
\mu_2\int_0^1z(x,1,t)\phi dx
{ }&{ }\leq\frac{\varepsilon}{c_p}\int_0^1\phi^2dx+\frac{\mu_2^2c_p}{4\varepsilon}\int_0^1z^2(x,1,t)dx \label{4:5a} \\
{ }&{ }\leq\varepsilon\int_0^1\phi_x^2dx+\frac{\mu_2^2c_p}{4\varepsilon}\int_0^1z^2(x,1,t)dx, \nonumber \\[1em]
\int_0^1f(\phi)\phi dx
\leq k_0\int_0^1|\phi|^\theta|\phi||\phi|dx
{ }&{ }\leq\|\phi\|_{2(\theta+1)}^\theta
\|\phi\|_{2(\theta+1)}\|\phi\|_2
\leq c_1\int_0^1\phi_x^2dx. \label{4:6a}
\end{align}
Combining \eqref{4:4}, \eqref{4:5a}, \eqref{4:6a} and \eqref{3:19} we obtain \eqref{4:3}.
\end{proof}

\begin{lemma}\label{L:4:3}
For any $\lambda_2,\tilde{\lambda}_2>0$, the functional $I_2$ satisfies,
\begin{equation}
\begin{split}
\frac{d}{dt}I_2(t)\leq (\mu_2\lambda_2-\delta)\int_0^1\phi_x^2dx+\bigg(J+\frac{\rho b^2}{4\mu^2\tilde{\lambda}_2}\bigg)\int_0^1\phi_t^2dx+\rho\tilde{\lambda}_2c_p\int_0^1u_t^2dx\\
+\frac{\mu_2^2c_p}{4\lambda_2}\int_0^1 z^2(x,1,t)dx-\int_0^1\hat{f}(\phi)dx.
\end{split}\label{4:7}
\end{equation}
\end{lemma}

\begin{proof}
Taking derivative of $I_2$, considering \eqref{2:5}, \eqref{2:6} and integrating by parts, we obtain
\begin{equation}
\begin{split}
\frac{d}{dt}I_2(t)=-\delta\int_0^1\phi_x^2dx+J\int_0^1\phi_t^2dx-\xi\int_0^1\phi^2dx+\rho\int_0^1 u_tw_tdx\\
+\mu\int_0^1w_x^2dx-\mu_2\int_0^1z(x,1,t)\phi dx-\int_0^1f(\phi)\phi dx.
\end{split}\label{4:8}
\end{equation}
By \eqref{4:6} we can get
\begin{eqnarray}
\mu\int_0^1w_x^2dx\leq\frac{b^2}{\mu}\int_0^1\phi^2dx\leq \xi c_p\int_0^1\phi_x^2dx,\label{4:9}\\
\int_0^1w_t^2dx\leq c_p\int_0^1w_{xt}^2dx\leq\frac{c_pb^2}{\mu^2}\int_0^1\phi_t^2dx\label{4:10}.
\end{eqnarray}
By using Young's and Poincaré's inequalities and \eqref{4:10}, we have
\begin{align}
\mu_2\int_0^1 z(x,1,t)\phi dx
&\leq \frac{\mu_2\lambda_2}{c_p}\int_0^1\phi^2dx+\frac{\mu_2c_p}{4\lambda_2}\int_0^1 z^2(x,1,t)dx \label{4:11} \\
&\leq \mu_2\lambda_2\int_0^1\phi_x^2dx +\frac{\mu_2c_p}{4\lambda_2}\int_0^1 z^2(x,1,t)dx, \nonumber \\[1ex]
\rho\int_0^1 u_tw_tdx
&\leq\rho\tilde{\lambda}_2c_p\int_0^1 u_t^2dx+\frac{\rho}{4\tilde{\lambda}_2c_p}\int_0^1 w_t^2dx \label{4:12} \\
&\leq\rho\tilde{\lambda}_2c_p\int_0^1 u_t^2dx+\frac{\rho b^2}{4\mu^2\tilde{\lambda}_2}\int_0^1\phi_t^2dx. \nonumber
\end{align}
From \eqref{4:8}-\eqref{4:12} we obtain \eqref{4:7}.
\end{proof}

\begin{lemma}\label{L:4:4}
Assume that \eqref{1:6} holds. Then the functional $I_3$ satisfies
\begin{multline}
\frac{d}{dt}I_3(t)\leq \delta\big[\phi_xu_x\big]_{x=0}^{x=1}-\frac{3\sqrt{\xi}}{8}\int_0^1\bigg(\frac{b}{\sqrt{\xi}}u_x+\sqrt{\xi}\phi\bigg)^2dx+
\bigg(J\sqrt{\xi}+\frac{\mu_1^2}{\sqrt{\xi}}\bigg)\int_0^1\phi_t^2dx\\
+\bigg(\frac{\xi^2c_p}{12\sqrt{\xi}}+\frac{6k_0^2}{\sqrt{\xi}}\bigg)\int_0^1\phi_x^2dx
+\frac{\mu_2^2}{\sqrt{\xi}}\int_0^1z^2(x,1,t)dx-\sqrt{\xi}\int_0^1\hat{f}(\phi)dx.
\label{4:13}
\end{multline}
\end{lemma}

\begin{proof}
Notice that
\begin{align}
\frac{d}{dt}I_3(t) =& J\int_0^1\!\!\!\phi_{tt}\bigg(\frac{b}{\sqrt{\xi}}u_x+\sqrt{\xi}\phi\bigg)dx+J\int_0^1\!\!\!\phi_t
\bigg(\frac{b}{\sqrt{\xi}}u_x+\sqrt{\xi}\phi\bigg)_tdx\nonumber \\
&\,\,\,+\frac{bJ}{\sqrt{\xi}}\int_0^1\!\!\!\phi_{xt} u_tdx +\frac{bJ}{\sqrt{\xi}}\int_0^1\!\!\!\phi_xu_{tt}dx.\label{4:14}
\end{align}
By using \eqref{2:5} and \eqref{2:6} and integration by parts, we have
\begin{multline}
\frac{d}{dt}I_3(t)=\delta[\phi_xu_x]_{x=0}^{x=1}-\sqrt{\xi}\int_0^1\bigg(\frac{b}{\sqrt{\xi}}u_x
+\sqrt{\xi}\phi\bigg)^2dx-\\ \mu_1\int_0^1\phi_t\bigg(\frac{b}{\sqrt{\xi}}u_x+\sqrt{\xi}\phi\bigg)dx
-\mu_2\int_0^1z(x,1,t)\bigg(\frac{b}{\sqrt{\xi}}u_x+\sqrt{\xi}\phi\bigg)dx\\-\frac{b}{\sqrt{\xi}}\int_0^1f(\phi)u_xdx
-\sqrt{\xi}\int_0^1f(\phi)\phi dx+J\sqrt{\xi}\int_0^1\phi_t^2dx.\label{4:15}
\end{multline}
By using Young's and Poincar\'e's inequality, we have

\begin{equation}\label{4:16}
\begin{array}{c}
-\mu_1\int_0^1\phi_t\bigg(\frac{b}{\sqrt{\xi}}u_x+\sqrt{\xi}\phi\bigg)dx\\ \leq\frac{\sqrt{\xi}}{4}
\int_0^1\bigg(\frac{b}{\sqrt{\xi}}u_x+\sqrt{\xi}\phi\bigg)^2dx
+\frac{\mu_1^2}{\sqrt{\xi}}\int_0^1\phi_t^2dx,
\end{array}
\end{equation}
\begin{equation}\label{4:17}
\begin{array}{c}
-\mu_2\int_0^1z(x,1,t)\bigg(\frac{b}{\sqrt{\xi}}u_x+\sqrt{\xi}\phi\bigg)dx
\\
\leq
\frac{\sqrt{\xi}}{4}
\int_0^1\bigg(\frac{b}{\sqrt{\xi}}u_x+\sqrt{\xi}\phi\bigg)^2dx
+\frac{\mu_2^2}{\sqrt{\xi}}\int_0^1z^2(x,1,t)dx,
\end{array}
\end{equation}
and by using \eqref{3:19}, for any $\varepsilon>0$
$$
\begin{array}{c}
-\frac{b}{\sqrt{\xi}}\int_0^1f(\phi)u_xdx\leq  \frac{k_0|b|}{\sqrt{\xi}}\|u_x\|_2\|\phi\|_{2(\theta+1)}^\theta\|\phi\|_{2(\theta+1)}\\

\leq\frac{k_0b^2\varepsilon}{\sqrt{\xi}}\int_0^1u_x^2dx+\frac{k_0}{4\sqrt{\xi}\varepsilon}\int_0^1\phi_x^2dx\\
\leq 3\varepsilon k_0\sqrt{\xi}\int_0^1\bigg(\frac{b}{\sqrt{\xi}}u_x+\sqrt{\xi}\phi\bigg)^2dx+
\bigg(\frac{2k_0\xi^2c_p\varepsilon}{\sqrt{\xi}}+\frac{k_0}{4\varepsilon\sqrt{\xi}}\bigg)\int_0^1\phi_x^2dx.
\end{array}$$
Taking $\varepsilon=1/24k_0$ in the previous inequality, we obtain
\begin{eqnarray}
-\frac{b}{\sqrt{\xi}}\int_0^1f(\phi)u_xdx\leq
\frac{\sqrt{\xi}}{8}\int_0^1\bigg(\frac{b}{\sqrt{\xi}}u_x+\sqrt{\xi}\phi\bigg)^2dx+
\bigg(\frac{\xi^2c_p}{12\sqrt{\xi}}+\frac{6k_0^2}{\sqrt{\xi}}\bigg)\int_0^1\phi_x^2dx.\label{4:18}
\end{eqnarray}
Combining \eqref{4:15}, \eqref{4:16}, \eqref{4:17} and \eqref{4:18} we obtain \eqref{4:13}.
\end{proof}

Next we deal with the boundary term in \eqref{4:13}. We define the function
\begin{eqnarray}
q(x)=2-4x, \quad x\in [0,1],\label{4:19}
\end{eqnarray}
so we have the following result.

\begin{lemma}\label{L:4:5}
For any $\varepsilon>0$, we have
\begin{multline}
\delta[u_x\phi_x]_{x=0}^{x=1}\leq-\frac{J\delta}{4\varepsilon}\frac{d}{dt}\int_0^1q\phi_t\phi_x dx-\frac{\rho\varepsilon}
{\mu}\frac{d}{dt}\int_0^1qu_tu_xdx+\frac{2\rho\varepsilon}{\mu}\int_0^1u_t^2dx\\
+\bigg(\frac{\delta^2}{2\varepsilon}+\frac{\delta^2}{4\varepsilon^3}+\frac{\delta^2}{2}+\frac{\varepsilon}
{\mu}+\frac{2\xi^2c_p\varepsilon}{\mu}+\frac{4\xi^2c_p\varepsilon}{b^2}\bigg)\int_0^1\phi_x^2dx
+\bigg(\frac{J\delta}{2\varepsilon}+\frac{\mu_1^2}{4\varepsilon^2}\bigg)\int_0^1\phi_t^2dx\\
+\bigg(\frac{\xi\varepsilon}{4}+\frac{3\xi\varepsilon}{\mu}+\frac{6\xi\varepsilon}{b^2}\bigg)\int_0^1\bigg(
\frac{b}{\sqrt{\xi}}u_x+\sqrt{\xi}\phi\bigg)^2dx+\frac{\mu_2^2}{4\varepsilon^2}\int_0^1z^2(x,1,t)dx.\label{4:20}
\end{multline}
\end{lemma}

\begin{proof}
By using Young's inequality we have
\begin{eqnarray}
\delta[u_x\phi_x]_{x=0}^{x=1}\leq \varepsilon[u_x^2(1)+u_x^2(0)]+\frac{\delta^2}{4\varepsilon}[\phi_x^2(1)+\phi_x^2(0)].\label{4:21}
\end{eqnarray}
Now, by using \eqref{2:6} and integrating by parts, we have
\begin{multline}
\frac{d}{dt}\int_0^1J\delta q\phi_t\phi_xdx=\\-\delta^2[\phi_x^2(1)+\phi_x^2(0)]+2\delta^2\int_0^1\phi_x^2dx
-\delta\sqrt{\xi}\int_0^1q\bigg(\frac{b}{\sqrt{\xi}}u_x+\sqrt{\xi}\phi\bigg)\phi_xdx\\
-\delta\mu_1\int_0^1q\phi_t\phi_xdx-\delta\mu_2\int_0^1qz(x,1,t)\phi_xdx
-J\delta\int_0^1qf(\phi)\phi_xdx+2J\delta\int_0^1\phi_t^2dx. \label{4:22}
\end{multline}
Using Young's inequality
\begin{gather}
\delta\sqrt{\xi}\int_0^1q\bigg(\frac{b}{\sqrt{\xi}}u_x+\sqrt{\xi}\phi\bigg)\phi_xdx\leq
\xi\varepsilon^2\int_0^1\bigg(\frac{b}{\sqrt{\xi}}u_x+\sqrt{\xi}\phi\bigg)^2dx
+\frac{\delta^2}{\varepsilon^2}\int_0^1\phi_x^2dx, \label{4:23} \\
\delta\mu_1\int_0^1q\phi_t\phi_xdx\leq\frac{\mu_1^2}{\varepsilon}\int_0^1\phi_t^2dx+\varepsilon
\delta^2\int_0^1\phi_x^2dx, \label{4:24} \\
\delta\mu_2\int_0^1qz(x,1,t)\phi_xdx\leq\frac{\mu_2^2}{\varepsilon}\int_0^1z^2(x,1,t)dx+\varepsilon
\delta^2\int_0^1\phi_x^2dx. \label{4:25}
\end{gather}
Combining \eqref{4:22}, \eqref{4:23}, \eqref{4:24}, \eqref{4:25} and \eqref{2:2} we have
\begin{multline}
\frac{d}{dt}\int_0^1J\delta q\phi_t\phi_xdx\leq-\delta^2[\phi_x^2(1)+\phi_x^2(0)]+\bigg(2\delta^2+
\frac{\delta^2}{\varepsilon^2}+2\varepsilon\delta^2\bigg)\int_0^1\phi_x^2dx\\
+\xi\varepsilon^2\int_0^1\bigg(\frac{b}{\sqrt{\xi}}u_x+\sqrt{\xi}\phi\bigg)^2dx
+\bigg(2J\delta+\frac{\mu_1^2}{\varepsilon}\bigg)\int_0^1\phi_t^2dx
+\frac{\mu_2^2}{\varepsilon}\int_0^1z^2(x,1,t)dx.\label{4:27}
\end{multline}
On the other hand, by using \eqref{2:5} and integration by parts, we have
\begin{eqnarray}
\frac{d}{dt}\int_0^1\rho qu_tu_xdx=-\mu[u_x^2(1)+u_x^2(0)]+2\mu\int_0^1u_x^2dx+b\int_0^1qu_x\phi_xdx
+2\rho\int_0^1u_t^2dx
\end{eqnarray}
and applying Young's inequality and \eqref{3:19} in the previous equality

\begin{multline}
\frac{d}{dt}\int_0^1\rho qu_tu_xdx\leq-\mu[u_x^2(1)+u_x^2(0)]+\bigg(3\xi+\frac{6\mu\xi}{b^2}\bigg)\int_0^1
\bigg(\frac{b}{\sqrt{\xi}}u_x+\sqrt{\xi}\phi\bigg)^2dx\\
+\bigg(1+2\xi^2c_p+\frac{4\mu\xi^2c_p}{b^2}\bigg)\int_0^1\phi_x^2dx+2\rho\int_0^1u_t^2dx.\label{4:28}
\end{multline}
Combining \eqref{4:27} and \eqref{4:28} we obtain \eqref{4:20}.
\end{proof}

Following the same line of argument given found in \cite{HOUARI-LASKRI} it is not difficult to obtain
\begin{lemma}\label{L:4:6}
The functional $I_4$ satisfies
\begin{eqnarray}
\frac{d}{dt}I_4(t)=-2I_4(t)-\frac{e^{2\tau}}{\tau}\int_0^1z^2(x,1,t)dx+\frac{1}{\tau}\int_0^1\phi_t^2dx
\label{4:30}.
\end{eqnarray}
\end{lemma}

Now we define the following Lyapunov functional $\mathcal{L}(t)$
\begin{equation}
\mathcal{L}(t):=ME(t)+\frac{\sqrt{\xi}}{8}\beta I_1(t)+NI_2(t)+I_3(t)+\frac{J\delta}{4\varepsilon}\int_0^1q\phi_t\phi_xdx
+\frac{\rho\varepsilon}{\mu}\int_0^1qu_tu_xdx+I_4(t)\label{4:31}
\end{equation}
where
\begin{eqnarray}
\beta:=\bigg[3\bigg(\frac{\mu\xi}{b^2}-1\bigg)+1\bigg]^{-1}.
\end{eqnarray}

\begin{lemma}\label{L:4:7}
Let $U(t)=(u(t),u_t(t),\phi(t),\phi_t(t),z(t))$ be the solution of problem \eqref{2:5}-\eqref{2:9}. For $M$ large enough, there exist two positives $\gamma_1$ and $\gamma_2$ depending on $M$, $N$ and $\varepsilon$ such that for any $t\geq 0$,
\begin{eqnarray}
\gamma_1E(t)\leq\mathcal{L}(t)\leq\gamma_2E(t).\label{4:32}
\end{eqnarray}
\end{lemma}

\begin{proof}
We have
\begin{align}
|\mathcal{L}(t)-ME(t)|
\leq{ }&{ }\frac{\sqrt{\xi}}{8}\beta \bigg| -\int_0^1(\rho uu_t+J\phi\phi_t)dx-\frac{\mu_1}{2}\int_0^1\phi^2dx\bigg| \label{4:33} \\
&{ }+N\bigg|\int_0^1(J\phi_t\phi+\rho u_tw)dx + \frac{\mu_1}{2}\int_0^1\phi^2dx\bigg| \nonumber \\
&{ }+\bigg|J\int_0^1\bigg(\frac{b}{\sqrt{\xi}}u_x + \sqrt{\xi}\phi\bigg)\phi_tdx + \frac{b J}{\sqrt{\xi}} \int_0^1\phi_x u_tdx\bigg| \nonumber \\
&{ }+\frac{J\delta}{4\varepsilon}\int_0^1|q\phi_t\phi_x|dx + \frac{\rho\varepsilon}{\mu}\int_0^1|qu_tu_x|dx \nonumber \\
&{ }+\bigg|\int_0^1\int_0^1e^{-2\tau y}z^2(x,y,t)dydx\bigg|. \nonumber
\end{align}
By using Young's and Poincaré's inequalities, \eqref{3:19} and \eqref{4:9} we obtain:
\newline
\resizebox{\textwidth}{!}{
\parbox{1.01\textwidth}{
\begin{align}
\frac{\sqrt{\xi}}{8}\beta \bigg| - \int_0^1(\rho uu_t+J\phi\phi_t)dx-\frac{\mu_1}{2}\int_0^1\phi^2dx\bigg|
\leq{ }&{ }\frac{\sqrt{\xi}\beta\rho}{16}\int_0^1u_t^2dx + \frac{\sqrt{\xi}\beta J}{16}\int_0^1\phi_t^2dx  \label{4:34} \\
&+\bigg[\frac{\sqrt{\xi}\beta c_p}{16}(J+\mu)+\frac{\xi^{5/2}\beta\rho c_p^2}{8b^2}\bigg]\int_0^1\phi_x^2dx \nonumber \\
&+ \frac{3\xi^{3/2}\beta\rho c_p}{16b^2}\int_0^1\bigg(\frac{b}{\sqrt{\xi}}u_x+\sqrt{\xi}\phi\bigg)^2dx, \nonumber \\[1em]
N\bigg|\int_0^1(J\phi_t\phi+\rho u_tw)dx+\frac{\mu_1}{2}\int_0^1\phi^2dx\bigg| 
\leq{ }&{ }\frac{N\rho}{2}\int_0^1u_t^2dx+\frac{NJ}{2}\int_0^1\phi_t^2dx \label{4:35} \\
&+\frac{Nc_p}{2}\bigg(J+\mu_1+\frac{\rho\xi c_p}{\mu}\bigg)\int_0^1\phi_x^2dx, \nonumber \\[1em]
\bigg|J\int_0^1\bigg(\frac{b}{\sqrt{\xi}}u_x+\sqrt{\xi}\phi\bigg)\phi_tdx+\frac{b J}{\sqrt{\xi}}\int_0^1\phi_x u_tdx\bigg|
\leq{ }&{ }\frac{Jb^2}{2}\int_0^1u_t^2dx +\frac{J}{2}\int_0^1\phi_t^2dx  \label{4:36} \\
&+\frac{J}{2}\int_0^1\bigg(\frac{b}{\sqrt{\xi}}u_x+\sqrt{\xi}\phi\bigg)^2dx \nonumber \\
&+\frac{J}{2\xi}\int_0^1\phi_x^2dx, \nonumber \\[1em]
\frac{J\delta}{4\varepsilon}\int_0^1|q\phi_t\phi_x|dx \leq{ }&{ }\frac{J\delta}{4\varepsilon}\int_0^1\phi_t^2dx
+\frac{J\delta}{4\varepsilon}\int_0^1\phi_x^2dx,\label{4:37} \\[1em]
\frac{\rho\varepsilon}{\mu}\int_0^1|qu_tu_x|dx 
\leq{ }&{ }\frac{\rho\varepsilon}{\mu}\int_0^1u_t^2dx + \frac{2\rho\xi^2 c_p\varepsilon}{\mu b^2}\int_0^1\phi_x^2dx \label{4:38A} \\
&+\frac{3\rho\xi\varepsilon}{\mu b^2}\int_0^1\bigg(\frac{b}{\sqrt{\xi}}u_x+\sqrt{\xi}\phi\bigg)^2dx. \nonumber
\end{align}
}}
From \eqref{4:33}-\eqref{4:38A}, we have
\begin{multline}
|\mathcal{L}(t)-ME(t)|\leq\alpha_1\frac{\rho}{2}\int_0^1u_t^2dx+\alpha_2\frac{J}{2}\int_0^1\phi_t^2dx+\alpha_3\frac{\delta}{2}
\int_0^1\phi_x^2dx\\
+\alpha_4\frac{1}{2}\int_0^1\bigg(\frac{b}{\sqrt{\xi}}u_x+\sqrt{\xi}\phi\bigg)^2dx
+\alpha_5\frac{\eta}{2}\int_0^1\int_0^1z^2(x,y,t)
dydx
\end{multline}
where
\begin{align}
\alpha_1&= \frac{\sqrt{\xi}\beta}{8}+N+\frac{Jb^2}{\rho}+\frac{2\varepsilon}{\mu},\\
\alpha_2&= \frac{\sqrt{\xi}\beta}{8}+N+1+\frac{\delta}{2\varepsilon},\\
\alpha_3&= \frac{\sqrt{\xi}\beta c_p(J+\mu)}{8\delta}+\frac{\xi^{5/2}\beta\rho c_p^2}{4\delta b^2}+\frac{Nc_p}{\delta}\bigg(J+\mu_1+\frac{\rho\xi c_p}{\mu}\bigg)+\frac{J}{\delta}+
\frac{J}{2\varepsilon}+\frac{4\rho\xi^2 c_p\varepsilon}{\mu\delta b^2},\\
\alpha_4&= \frac{3\xi^{3/2}\beta\rho c_p}{8b^2}+J+\frac{6\rho\xi\varepsilon}{\mu b^2},\\
\alpha_5&= \frac{2}{\eta}.
\end{align}
By considering
\begin{eqnarray}
\tilde{C}=\max\{\alpha_i;\ i=1,\ldots,5\}
\end{eqnarray}
we arrive
\begin{eqnarray}
|\mathcal{L}(t)-ME(t)|\leq \tilde{C}E(t).
\end{eqnarray}
So, we can choose $M$ large enough so that $\gamma_1=M-\tilde{C}>0$ and $\gamma_2=M+\tilde{C}>0$ such that \eqref{4:32} holds.
\end{proof}

We are ready to prove the main result of this paper.

\begin{proof}[Proof of Theorem \ref{T:4:1}]
It follows from \eqref{3:20}, \eqref{4:3}, \eqref{4:7}, \eqref{4:13}, \eqref{4:20}, and \eqref{4:30}
\begin{align}
\frac{d}{dt}\mathcal{L}(t)&\leq \bigg[-\frac{\rho\sqrt{\xi}\beta}{8}+N\rho\tilde{\lambda}_2c_p+
\frac{2\rho\varepsilon}{\mu}\bigg]\int_0^1u_t^2dx\nonumber\\
& +\bigg[\!-\!MC_E\!-\!\frac{\sqrt{\xi}\beta J}{8}\!+\!N\bigg(J+\frac{\rho b^2}{4\mu^2\tilde{\lambda}_2}\bigg)
\!+\!\bigg(J\sqrt{\xi}\!+\!\frac{\mu_1^2}{\sqrt{\xi}}\bigg)\!+\!\bigg(\frac{J\delta}{2\varepsilon}\!+\!\frac{\mu_1^2}{4\varepsilon^2}
\bigg)\!+\!\frac{1}{\tau}\bigg]\int_0^1\!\!\phi_t^2dx\nonumber\\
& + \frac{\sqrt{\xi}\beta}{8}\bigg[\delta+\varepsilon+c_1+\frac{2\xi^2c_p}{b^2}
\bigg(\mu-\frac{b^2}{\xi}\bigg)\bigg] \int_0^1\phi_x^2dx \nonumber\\
\end{align}
\begin{align}
& +\bigg[N(\mu_2
\lambda_2-\delta)+\frac{\xi^2c_p}{12\sqrt{\xi}}+\frac{6k_0}{\sqrt{\xi}}+\frac{\delta^2}{2\varepsilon}
+\frac{\delta^2}{4\varepsilon^2}+\frac{\delta^2}{2}+\frac{\xi}{\mu}+\frac{2\xi^2c_p\varepsilon}{\mu}
+\frac{4\xi^2c_p\varepsilon}{b^2}\bigg] \int_0^1\phi_x^2dx\nonumber\\
& +\bigg[-\frac{\sqrt{\xi}}{4}+\frac{\xi\varepsilon}{4}+
\frac{3\xi\varepsilon}{\mu}+\frac{6\xi\varepsilon}{b^2}\bigg]\int_0^1\bigg(\frac{b}{\sqrt{\xi}}u_x+
\sqrt{\xi}\phi\bigg)^2dx \nonumber\\
& +\bigg[-MC_E+\frac{\mu_2^2c_p}{32\varepsilon}+\frac{N\mu_2^2c_p}{4\lambda_2}+\frac{\mu_2^2}{\sqrt{\xi}}
+\frac{\mu_2^2}{4\varepsilon^2}-\frac{e^{2\tau}}{\tau}\bigg]\int_0^1z^2(x,1,t)dx \nonumber\\&-(N+\sqrt{\xi})
\int_0^1\hat{f}(\phi)dx -e^{-2\tau}\int_0^1\int_0^1z^2(x,y,t)dydx.\label{4:40}
\end{align}
At this point, we have to choose our constants very carefully. First, let us choose $\lambda_2$ small enough such that
\begin{eqnarray}
\lambda_2<\frac{\delta}{\mu_2}.
\end{eqnarray}
Then, take $\varepsilon$ small enough such that
\begin{eqnarray}
\varepsilon\leq\min\bigg\{\frac{\sqrt{\xi}}{8}\bigg/\bigg(\frac{\xi}{4}+\frac{3\xi}{\mu}+\frac{6\xi}{b^2}\bigg),
\frac{\sqrt{\xi}\mu\beta}{64}\bigg\}.
\end{eqnarray}
Then, we select $N$ large enough so that
\begin{align*}
N(\delta-\mu_2\lambda_2) >& 
\frac{\sqrt{\xi}\beta}{8}\bigg[\delta+\varepsilon+c_1+\frac{2\xi^2c_p}{b^2}
\bigg(\mu-\frac{b^2}{\xi}\bigg)\bigg]\\
&\,\,+\frac{\xi^2c_p}{12\sqrt{\xi}}+\frac{6k_0}{\sqrt{\xi}}+\frac{\delta^2}{2\varepsilon}
+\frac{\delta^2}{4\varepsilon^2}+\frac{\delta^2}{2}+\frac{\xi}{\mu}+\frac{2\xi^2c_p\varepsilon}{\mu}
 +\frac{4\xi^2c_p\varepsilon}{b^2}.
\end{align*}
Now, we pick $\tilde{\lambda}_2$ so small that
\begin{eqnarray*}
\tilde{\lambda}_2<\frac{\sqrt{\xi}\beta}{32 Nc_p}.
\end{eqnarray*}
Finally, we choose $M$ large enough so that, there exists a positive constant $\beta_1>0$, such that \eqref{4:40} becomes
\begin{align}
\frac{d}{dt}\mathcal{L}(t)\leq& -\beta_1\int_0^1\bigg[u_t^2+\phi_t^2+\phi_x^2+\bigg(\frac{b}{\sqrt{\xi}}u_x+\sqrt{\xi}\phi
\bigg)^2+z^2(x,1,t)+\hat{f}(\phi)\bigg]dx\nonumber\\
& -\beta_1\int_0^1\int_0^1z^2(x,y,t)dydx\label{4:41}
\end{align}
which implies that there exists also $\beta_2>0$, such that
\begin{eqnarray}
\frac{d}{dt}\mathcal{L}(t)\leq-\beta_2E(t),\quad\forall t\geq 0,
\end{eqnarray}
and from \eqref{4:32}
\begin{eqnarray}
\frac{d}{dt}\mathcal{L}(t)\leq-\frac{\beta_2}{\gamma_1}\mathcal{L}(t),\quad\forall t\geq 0.\label{4:42}
\end{eqnarray}
Integrating \eqref{4:43} and using again \eqref{4:32} we arrive
\begin{eqnarray}
E(t)\leq-\frac{\gamma_2}{\gamma_1}E(0)e^{-(\beta_2/\gamma_2)t},\quad\forall t\geq 0,\label{4:43}
\end{eqnarray}
which gives us that the exponential stability holds for any $U_0\in D(A)$. Noting that $D(A)$ is
dense $\mathcal{H}$, we can extend the energy inequalities to phase space $\mathcal{H}$. Thus we complete
the proof of the theorem.
\end{proof}


\EditInfo{June 07,2019}{April 13, 2022}{Durvudkhan Suragan}

\end{paper}